\documentclass[12pt]{article}
\usepackage{latexsym}
\usepackage{graphicx}

\usepackage{amsmath, amsthm, amsfonts, amssymb}

\title{Residual Velocities in Steady Free Boundary Value Problems
of Vector Laplacian Type}

\author{Wan Chen \thanks{Department of Mathematics, UBC,
wanchen@math.ubc.ca} \and Brian Wetton \thanks{corresponding author,
Department of Mathematics, UBC, wetton@math.ubc.ca.}}

\begin{document}
\maketitle

\abstract{This paper describes a technique to determine the linear
well-posedness of a general class of vector elliptic problems that
include a steady interface, to be determined as part of the problem,
that separates two subdomains. The interface satisfies mixed Dirichlet
and Neumann conditions. We consider ``2+2'' models, meaning two
independent variables respectively on each subdomain. The governing
equations are taken to be vector Laplacian, to be able to make
analytic progress. The interface conditions can be classified into
four large categories, and we concentrate on the one with most
physical interest. The well-posedness criteria in this case are
particularly clear. In many physical cases, the movement of the
interface in time-dependent situations can be reduced to a normal
motion proportional to the residual in one of the steady state
interface conditions (the elliptic interior problems and the other
interface conditions are satisfied at each time). If only the steady
state is of interest, one can consider using other residuals for the
normal velocity. Our analysis can be extended to give insight into
choosing residual velocities that have superior numerical
properties. Hence, in the second part, we discuss an iterative method
to solve free boundary problems. The advantages of the correctly
chosen, non-physical residual velocities are demonstrated in a
numerical example, based on a simplified model of two-phase flow with
phase change in porous media.}

\vspace{5mm}


\section{Introduction}

Free boundary problems (FBPs) have motivated several studies in the
past due to their vast applications in fluid flow, phase change
models and other fields. Some classical free boundary problems are:
the dam seepage problem \cite{crank}, incompressible two-phase 
flow \cite{droplet1,droplet2} ({\em i.e.} falling droplets 
or rising
bubbles \cite{bubble}, the Alt-Caffarelli problem \cite{alt}, the
classical Stefan problem \cite{crank}, etc.  From a mathematical
point of view, FBPs are boundary value problems with an unknown
boundary. The motion (unsteady case) or position (steady case) of
the boundary has to be determined together with the solution of the
given partial differential equations on one (free surface) or both
sides (interface) of the free boundary. The coupling of the free
boundary to the interior is always nonlinear \cite{shape1}, and thus
FBPs are often not easy to solve.

The common structure of all the examples above is that at steady
state, they all have second order elliptic governing equations for $m$
unknowns on one side and $n$ unknowns on the other side of the free
boundary. We denote this situation as an ``m+n'' problem. There must
then be $m+n+1$ Dirichlet-Neumann conditions at the interface.  For
fourth order problems such as biharmonic equation, we would need
$2m+2n+1$ conditions at the interface to determine an ``m+n''
problem. By this generalization, many more complicated problems can be
formulated. However, fourth order problems and second order elliptic
systems other than vector Laplacian are not considered further in this
present study.

To solve free boundary problems numerically there are three main
kinds of methods: capturing methods, front tracking methods and
level set methods. Capturing methods are based on Eulerian
formulation and the problems are reformulated and solved in the
whole domain.  The interface location is recovered from the discrete
solution. In these methods, the interface conditions are not
specified explicitly. A classical example is the Enthalpy method
\cite{crank}. The alternative is to discretize the interface
explicitly \cite{tracking1} \cite{tracking} or via a level set
approach \cite{level}. In both these cases, the interface conditions
(and interface velocity for time-dependent problems) are implemented
explicitly. This can be done by considering the domains as disjoint
and discretizing the equations and interface conditions directly
\cite{direct} or by discretizing the entire domain and modifying the
discretization near the interface \cite{modified}. The latter
approach combined with modern level set techniques can also be
considered a capturing method.

In a few cases, steady state solutions can be obtained directly
using shape optimization \cite{shape1,shape2}. Another approach to reach
steady state solutions is to solve the transient time-dependent
problem to long time. In the class of problems we consider, the
normal motion of the interface is driven by the difference of
solutions on either side of the interface. For example,
solidification boundaries are driven by a Stefan condition
\cite{crank}. It should be made clear that we are not considering 
here the problem of dendritic growth \cite{al,caflisch} which in our 
framework would be ill-posed, but is regularized by higher order 
(curvature) terms. We consider only well-posed problems which do 
not need regularization. 
In the physical example of section \ref{sphys}, the
interface is moved by the net mass flux at the interface. Those are
typical conditions that give normal velocity of each point on the
interface. In these examples and others of physical interest, the
normal interface velocity is given by the residual in one of the
steady state interface conditions. In this paper, we apply the
residual velocity method (a tracking technique) of \cite{roger}:
\begin{enumerate}
\item Given an initial interface $\Gamma$ and tracking points
$\{\mathbf{x}: \mathbf{x} \in \Gamma\}$, solve numerically a fixed
boundary value problem satisfying all of the steady interface
conditions but one.
\item Substitute the solutions into the unsatisfied condition and
find the residual $R(\mathbf{x})$.
\item If the residual is larger than a tolerance, explicitly evolve
the interface using the residual as a normal velocity of the tracking
points:
\[
\mathbf{V}_n = R(\mathbf{x})
\]
\item Repeat the process until the residual is less than a given tolerance.
\end{enumerate}
At termination, all interface conditions are approximately
satisfied. This is an example of a {\em value} method \cite{value}
(since it uses only values of the solution not shape derivatives to
update the interface). In this algorithm, the choice of residual
velocity is not specified, so we can choose different interface
conditions and their combinations as the residual velocity, not only
the real, physical velocity. As a result, we can show both analytically
and numerically that the residual velocity chosen by our criteria has
better numerical properties than the physical one. Specifically, our
method can use time steps independent of mesh size in the explicit
step 3 above. Because of the use of an artificial interface velocity,
our method provides accurate solutions only at steady state.

There are few publications addressing the analytic theory of the
general class of the steady free boundary problems considered here. In
\cite{roger} the authors present a linear theory of the well-posedness
of the general class of ``1+1'' Laplacian problems. Our work is an
extension of \cite{roger}. Specifically, the ideas considered there
for the scalar case are extended to the ``2+2'' case. An analytic
simplification to the well-posedness is identified in a physically
important class of problems. The class is represented by a simplified
model of two-phase flow with phase change in porous media. In
addition, a different numerical implementation of the corresponding
residual velocity method is given. While the numerical method is not
as general as that of \cite{roger} it does show that accurate results
can be obtained without grid refinement at the interface. Attaining
this in more general codes is a goal for future work.

The framework of this paper is following: In section \ref{sclass}, we
discuss the four classes of second order ``2+2'' problems with
different combinations of Dirichlet-Neumann conditions. A physically
interesting example, a simplified model of two-phase flow with phase
change in porous media, is presented in section \ref{sphys}. In
section \ref{sposed}, well-posedness criteria of ``2+2'' Laplacian
equations from \cite{roger} are reviewed and applied to the example
and generalized. The extension to 3-D problem has been shown as
well. Then in section \ref{snum}, the example of vector Laplacian type
free boundary problem are computed numerically, using a finite
difference method in mapped coordinates (a so-called fictitious 
domain method \cite{shape2}). The performances of both
physical velocity and our residual velocity are compared, and the
results agree with the analytical predictions.


\section{Class Division of ``2+2'' Problems}
\label{sclass}

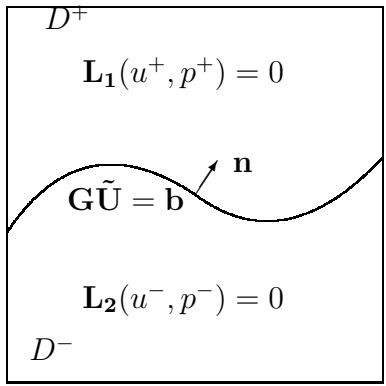
\begin{figure}
\setlength{\unitlength}{1cm} \hspace{4.5cm}
\begin{picture}(0,3.5)
\put(0,0){\line(0,1){5}} \put(0,0){\line(1,0){5}}
\put(5,0){\line(0,1){5}} \put(0,5){\line(1,0){5}}
\qbezier(0,2)(1,3.5)(2.5, 2.5)
\qbezier(2.5, 2.5)(3.7,1.6)(5,3)
\put(2.5,2.5){\vector(2,3){0.3}} \put(3,2.8){$\mathbf{n}$}
\put(0.5,4.7){$D^+$}
\put(0.3,0.3){$D^-$}
\put(1,4){$\mathbf{L_1} (u^+,p^+) =0$}
\put(1,1){$\mathbf{L_2} (u^-,p^-) =0$}
\put(0.8, 2.3) {$\mathbf{G\tilde{U}}=\mathbf{b}$}
\end{picture}
\caption{\label{fig1} The model of second order ``2+2'' problems.
$\mathbf{L_{1,2}}$ are linear, second order elliptic operators and
there are five interface conditions represented by $\mathbf{G}$.}
\end{figure}

Consider a two-dimensional model, when $m=n=2$. We have a ``2+2''
vector elliptic problems for variables $u^{\pm},p^{\pm}$, assuming
that $\mathbf{L_1} (u^+,p^+)=0$ in the upper subdomain $D^+$ and
$\mathbf{L_2} (u^-,p^-)=0$ in the lower subdomain $D^-$, as shown in
Fig.~\ref{fig1}. There are five conditions at the interface. Here,
we consider $\mathbf{L_1}$ and $\mathbf{L_2}$ to be general, vector,
linear, second order elliptic operator but in our analysis below we
consider only the case of $\mathbf{L_{1,2}}$ both being the vector
Laplacian.

The five Dirichlet-Neumann boundary conditions at $\Gamma$ in matrix
form, will be:
\begin{equation}
\label{gub}
\mathbf{G\tilde{U}}\;=\;\mathbf{b}.
\end{equation}
where $\mathbf{G}$ is a $5\times 8$ matrix, and
$\mathbf{\tilde{U}^T}=(u^+_\mathbf{n},\;\;u^-_\mathbf{n},
\;\;p^+_\mathbf{n},\;\;p^-_\mathbf{n},\;\;u^+,\;\;u^-,\;\;p^+,\;\;p^-)$
is a vector composed of all unknown variables and their normal
($\mathbf{n}$ pointing from $D^-$ to $D^+$) derivatives at the
interface. The matrix $\mathbf{G}$ can be split naturally into two
parts:
\[
\mathbf{G}=(G_N\;\;|\;\;G_D )
\]
where $G_N$ and $G_D$ are $5\times 4$ sub-matrices, representing
respectively Neumann conditions and Dirichlet conditions.

We use the rank $r$ of $G_N$ (or equivalently the number $5-r$ of pure
Dirichlet conditions) as the criteria to classify possible
forms of $\mathbf{G}$. As a result, there are four classes of boundary
conditions corresponding to $r$=4, 3, 2, 1 below. Normal forms can be
obtained using the following allowable operations:
\begin{enumerate}
\item Row operations on $\mathbf{G}$.
\item Simultaneous swap of columns 1 and 2, 3 and 4, 5 and 6, 7 and 8. 
This swap is accompanied by a sign change of columns 1-4. 
(corresponds to a relabelling of $D^-$ and $D^+$).
\item Column operations between columns 1 and 3 with identical
operations on columns 5 and 7 (linear combination of unknowns in
$D^+$).
\item Column operations between columns 2 and 4 with identical
operations on columns 6 and 8 (linear combination of unknowns in
$D^-$).
\end{enumerate}
The following generic forms result:
\begin{description}
\item[Class A:] One pure Dirichlet condition, and thus the rank of
$G_N$ is $4$:
\[
\mathbf{G}=\left(
\begin{array} {lll}
 1\;\;0\;\;0\;\;0\;\;|\;\;0\;\; & G_{1,6}\;& \;G_{1,7}\;\;G_{1,8} \\
 0\;\;1\;\;0\;\;0\;\;|\;\;0\;\; & G_{2,6}\;& \;G_{2,7}\;\;G_{2,8} \\
 0\;\;0\;\;1\;\;0\;\;|\;\;0 \;\;& \cdots\;& \;\cdots\;\;\cdots    \\
 0\;\;0\;\;0\;\;1\;\;|\;\;0\;\; & \cdots\;& \;\cdots\;\;\cdots    \\
 0\;\;0\;\;0\;\;0\;\;|\;\;1\;\; & 0 \;& \;G_{5,7}\;\;G_{5,8}
\end{array}
\right)
\]
\item[Class B:] Two pure Dirichlet conditions, and the rank of
$G_N$ is $3$.
\[
\mathbf{G}=\left(
\begin{array} {lllll}
 1\;0\;\;0\;\;& G_{1,4}\;& \;|\;\;0\;& \;0\;& \; G_{1,7}\;\;G_{1,8} \\
 0\;1\;\;0\;\;& G_{2,4}\;& \;|\;\;0\;& \;0\;& \; G_{2,7}\;\;G_{2,8} \\
 0\;0\;\;1\;\;& G_{3,4}\;& \;|\;\;0\;& \;0\;& \; G_{3,7}\;\;G_{3,8} \\
 0\;\;0\;\;0\;& \;0    \;& \;|\;\;G_{4,5}& G_{4,6} & G_{4,7}\;\;G_{4,8} \\
 0\;\;0\;\;0\;& \;0    \;& \;|\;\;G_{5,5}& G_{5,6} & G_{5,7}\;\;G_{5,8}
\end{array}
\right)
\]
where the last two rows are linearly independent but no further
reduction is possible.
\item[Class C:] Three pure Dirichlet conditions:
\[
\mathbf{G}=\left(
\begin{array} {ccccl}
 G_{1,1}& G_{1,2} & G_{1,3} & G_{1,4} & \;|\;\;0\;\;0\;\;0\;\;G_{1,8} \\
 G_{2,1}& G_{2,2} & G_{2,3} & G_{2,4} & \;|\;\;0\;\;0\;\;0\;\;G_{2,8} \\
 0\;\;  & 0       &0 & 0\;& \;|\;\;1\;\;0\;\;0\;\;G_{3,8} \\
 0\;\;  & 0       &0 & 0\;& \;|\;\;0\;\;1\;\;0\;\;G_{4,8} \\
 0\;\;  & 0       &0 & 0\;& \;|\;\;0\;\;0\;\;1\;\;G_{5,8}
\end{array}
\right)
\]
where the first two row vectors restricted to the left block are
linearly independent but no further reduction is possible.
\item[Class D:] Four pure Dirichlet conditions:
\[
\mathbf{G}=\left(
\begin{array} {cccc}
 1\;\;  0\;\;  &G_{1,3}&G_{1,4} & \;|\;\;0\;\;0\;\;0\;\;0  \\
 0\;\;  0\;\;  &0\;\;  &0\;& \;|\;\;1\;\;0\;\;0\;\;0 \\
 0\;\;  0\;\;  &0\;\;  &0\;& \;|\;\;0\;\;1\;\;0\;\;0 \\
 0\;\;  0\;\;  &0\;\;  &0\;& \;|\;\;0\;\;0\;\;1\;\;0 \\
 0\;\;  0\;\;  &0\;\;  &0\;& \;|\;\;0\;\;0\;\;0\;\;1
\end{array}
\right)
\]
\end{description}

Class C is arguably of most physical interest. If the vector elliptic
equations represent conservation of two quantities with fluxes
proportional to gradients, then the conservation of these quantities
at the interface are represented by exactly two pure ($G_{1,8}$ and
$G_{2,8}$ both zero) Neumann conditions. We define the subclass
$\tilde{C}$ to be the one with these additional conditions. In the
following section, we present a physical example in this class
originally derived from \cite{lloyd}. Class D always leads to a single
pure Neumann condition. We define subclasses $\tilde{A}$ and
$\tilde{B}$ to be the ones with four and three, respectively, pure
Neumann conditions.


\section{Physical example}
\label{sphys}

We consider now a physical problem involving two-phase flow with phase
change in porous media. Consider a closed system consisting of a sand
pack, water vapour and liquid water. If the system is heated at the
top, vapour is formed, and in some cases, two distinct zones with free
boundary will be observed. In the upper zone $D^+$, there is only
vapour, while in the lower zone $D^-$, two-phase vapour and liquid water
can be found. This model is studied in \cite{lloyd}.

In the upper region $D^+$, the variable $P_v(x,t)$ and temperature
$T^+(x,t)$ describe the vapour region, and the governing equations
are:
\begin{eqnarray}
\label{eq22}
\Delta T^+ & = & 0 \\
\label{eq23}
\nabla \cdot (\rho_v \mathbf{u}_v) & = & 0
\end{eqnarray}
where $\rho_v$ and $\mathbf{u}_v$ denote density and velocity of the
vapour respectively, which are the functions of temperature $T^+$ and
pressure $P_v$. According to Darcy's law, we have:
\[
\mathbf{u}_v=-\frac{\kappa}{\mu_v}\nabla P_v
\]
with $\kappa$ the permeability and $\mu_v$ the dynamic viscosity. If
we assume the gas is ideal, then:
\[
\rho_v=\frac{M}{R}\frac{P_v}{T^+}
\]
where $M$ is molar mass of water, $R$ is the universal gas constant.
So the equation (\ref{eq23}) becomes:
\begin{equation}
\label{eq24}
\nabla \cdot(\frac{P_v}{T^+} \nabla P_v)=0
\end{equation}
Then (\ref{eq22}) and (\ref{eq24}) give an elliptic problem for $P_v$
and $T^+$.

Flow in the lower two-phase region is described by temperature $T^-$
and saturation $s(x,y)$. The saturation is the fraction of pore
space occupied by liquid water. Energy conservation and mass
conservation give us:
\begin{eqnarray}
\label{eq25}
\nabla \cdot(K\nabla T^-)-
h_{vap}\nabla \cdot(\rho_v \mathbf{u}_v) & = & 0 \\
\label{eq26}
\nabla \cdot (\rho_l \mathbf{u}_l+\rho_v \mathbf{u}_v) & = & 0
\end{eqnarray}
where $\mathbf{u}_l,\mathbf{u}_v$ and $\rho_v$ are functions of
$T^-$ and $s$. Darcy's law includes the relative permeabilities as
well as the capillary pressure, which we assume are functions of the
$s$:
\[
\mathbf{u}_v=-\frac{\kappa}{\mu_v}(1-s)^3\nabla P_{sat}(T^-)
\]
\[
\mathbf{u}_l=-\frac{\kappa}{\mu_l} s^3 \nabla
(P_{sat}(T^-)-P_c(s))
\]
where $\mu_v, \mu_l$ are respectively dynamical viscosities of the
vapour and liquid, $P_{sat}$ is the saturation vapour pressure, and $P_c
= P_v - P_l$ is the capillary pressure. Since evaporation and
condensation rates scale to be very large in a porous media, we assume
$P_v = P_{sat}(T^-)$. Therefore in the two-phase region, we have an
elliptic problem for $T^-$ and $s$.

On the free boundary, There are five Dirichlet-Neumann conditions
required for the problem of four variables:
\[
\begin{aligned}
 s &= 0\;\;\;& \mbox{saturation is zero at the interface} \\
 [T]&=0\;\;\;& \mbox{temperature is continuous} \\
 [P_v]&=0\;\;\;&\mbox{vapour pressure is continuous } \\
 (\rho_v \mathbf{u}_v)^+\cdot {\mathbf{n}} &=(\rho_v \mathbf{u}_v+\rho_l \mathbf{u}_l)^- \cdot
 {\mathbf{n}}& \mbox{mass conservation} \\
 [K T_\mathbf{n}]&= h_{vap}(\rho_l \mathbf{u}_l)\cdot{\mathbf{n}}\;\;&\mbox{heat conservation}
 \end{aligned}
 \]
where $[\cdot]$ denotes the difference between counterparts on each
side of the interface. For example: $[T]=T^+ - T^-$.

In order to make progress on this problem in our framework,
considerable simplification is needed.  Specifically, we need to
simplify the problem to get a vector Laplacian operator on each
subdomain and linear interface conditions. The model derived below
lacks several aspects of physical interest and mathematical
difficulty. However, it does retain the underlying structure of two
elliptic systems coupled by five mixed Dirichlet-Neumann conditions at
the interface. First, we set $P_v/T^+=1$ in (\ref{eq24}), 
which means the
vapour density doesn't change. This leads to Laplacian equation of
pressure: $\Delta P=0$, then in the upper region $D^+$, we have:
\[
\Delta P=0, \mbox{\ \ \ } \Delta T^+=0
\]

In the lower region $D^-$, we simplify using the assumption that
evaporation or condensation only happens at the interface.
Moreover, we assume the relative permeability of vapour and liquid in
$D^-$ is independent of saturation $s$, the capillary pressure $P_c$
is linear in $s$ and saturation vapour pressure $P_{sat}$ is linear in
$T^-$. In the end, (\ref{eq25}) and (\ref{eq26}) can be reduced to two
Laplacian equations:
\[
\Delta s=0, \mbox{\ \ \ } \Delta T^-=0.
\]

Correspondingly, the interface conditions are given:
\[
\begin{aligned}
 s &= 0\;\;\;& \mbox{saturation is zero across interface} \\
 [T]&=0\;\;\;& \mbox{temperature is continuous} \\
 P&=T^-\;\;\;&\mbox{Vapour pressure is continuous } \\
 [K T_\mathbf{n}]&= -s_\mathbf{n} \;\;\;&\mbox{heat flux is continuous } \\
 P_\mathbf{n} &= T^-_\mathbf{n}+ s_\mathbf{n}
 \;\;\;&\mbox{mass conservation}
 \end{aligned}
\]
Note that in the example, the interface is featured by $s = 0$, and
its movement is driven by the mass flux going through the interface,
so we say the physical residual velocity comes from the last Neumann
condition. The five conditions can be written in matrix form (\ref{gub})
where $\mathbf{G}$ and $\mathbf{\tilde{U}}$ are respectively:
\begin{equation}
\label{eq28}
\mathbf{G}= \left(
\begin{array} {cccccccc}
 0  &0  &0  &0\;\;&\;\;0 &0  &0 &1        \\
 0  &0  &0  &0\;\;&\;\;1 &-1 &0 &0       \\
 0  &0  &0  &0\;\;&\;\;0 &1 &-1 &0       \\
 K^+  &-K^-  &0  &+1\;\;&\;\;0 &0 &0 &0     \\
 0  &-1  &1  &-1\;\;&\;\;0 &0 &0 &0
\end{array}
 \right)
\end{equation}
and
\[
\mathbf{\tilde{U}^T}=(
T^+_\mathbf{n},\;T^-_\mathbf{n},\;P_\mathbf{n},\;s_\mathbf{n},
\;T^+,\;T^-,\;P,\;s)
\]
Clearly, the boundary conditions are composed of three pure
Dirichlet and two pure Neumann conditions, which belongs to Class
$\tilde{C}$ identified in the last section.


\section{Linear Well-Posedness and Residual Velocities}
\label{sposed}

In the first subsection below we review briefly the theory developed
in \cite{roger} targeted to the ``2+2'' case. Later, we identify
situations where the analysis gives particularly clear insight. These
include the physically important case $\tilde{C}$ identified in the
previous section. Application of the theory to the example above is
done, and a residual velocity with superior numerical properties is
identified.

\subsection{General Theory}

\subsubsection{Base Solutions}
\label{sbase}

We find base solutions in 2D (coordinates $x$ and $y$) with a flat
interface $y=0$. Take $\mathbf{u}_0^T = (u^+, u^-, p^+, p^-)$ to be
the base solution of the FBP satisfying:
\begin{equation}
\label{eq1}
\Delta \mathbf{u}_0=0 \qquad\mbox{in $y>0$ and $y<0$}
\end{equation}
with interface conditions:
\begin{equation}
\label{eq2}
\mathbf{G\,\tilde{U}_0}=\mathbf{b} \qquad\mbox{on the
interface $y=0$}
\end{equation}
where $\mathbf{\tilde{U}}_0=\left(\begin{array}{c}
\mathbf{u}_{0\mathbf{n}} \\ \mathbf{u}_0 \end{array}\right)$.

We consider base solutions independent of $x$ with $\mathbf{b}$
constant, so (\ref{eq1}) can be reduced to $\mathbf{u}_{0yy}=0$ and
$\mathbf{u}_0$ is linearly dependent on y. That yields:
\[
\mathbf{u}_0=y\;\mathbf{r}^1+\mathbf{r}^0
\]
where
\[
\mathbf{\tilde{U}_0}=\left(\begin{array}{c} \mathbf{r}^{1}
\\ \mathbf{r}^{0} \end{array}\right)
\]
is a constant vector solving (\ref{eq2}). The base solution is
composed of a particular solution and homogeneous solution
\[
\mathbf{\tilde{U}_0}=\mathbf{\tilde{r}}_H+\mathbf{\tilde{r}}_P,
\]
with $\mathbf{G\tilde{r}}_H=\mathbf{0}$. Hence, $\mathbf{\tilde{r}}_H$
belongs to the nullspace of $\mathbf{G}$, and
$\mathcal{N}(\mathbf{G})$ has rank 3.

In order to specify a unique solution, we need three more conditions
denoted by a vector $\mathbf{q}=(q_1,q_2,q_3)^T$, corresponding to the
ordering of a basis of vectors of $\mathcal{N}(\mathbf{G})$, so
\[
\mathbf{\tilde{U}}_0 =
\left(\begin{array}{c}
\mathbf{r}^{1} \\
\mathbf{r}^{0}
\end{array}\right) = \mathcal{N}(\mathbf{G})\cdot \mathbf{q}
\]
where we mean (with a slight abuse of notation) by
$\mathcal{N}(\mathbf{G})$ the $8 \times 3$ matrix with a basis for
this space in columns.

Consider the base solution found above to describe the local behaviour
of the interface near a point, where the coordinate system has been
rotated so that the local tangent plane is $y=0$. This local base
solution depends on the local data $\mathbf{b}$ and also the vector
$\mathbf{q}$, which we refer to as the global fluxes. This
represents the dependence of the local solution on far field
conditions.


\subsubsection{Linear well-posedness}

We consider how small perturbation of boundary conditions about a flat
interface affect the solution:
\begin{equation}
\label{eq34}
\mathbf{G\tilde{U}}=\mathbf{b}+\epsilon \mathbf{f} e^{i\alpha x}
\end{equation}
where $\mathbf{f}$ is a constant, $\alpha$ is the Fourier mode, and
$\epsilon$ is a small parameter. As a result, we will obtain a new
solution and new free boundary. If we can show that the linearized
solution depends continuously on the data $\mathbf{f}$ for each
$\alpha$ we say the problem is linearly well-posed.

After we perturb the interface condition using (\ref{eq34}), we expect
the new solution $\mathbf{u}$ and free boundary $y = \eta (x)$ to
be in the form:
\[
\begin{aligned}
\mathbf{u}(x,y) &= \mathbf{u}_0 (y) + \epsilon \mathbf{u}_1
(x,y)+ O(\epsilon^2)\\
\eta(x) &= \epsilon \eta_1(x)+ O(\epsilon^2)
\end{aligned}
\]
where $\mathbf{u}_0, \eta_0(x)$ are the base solution from section
\ref{sbase} above. Following \cite{roger} it is known that linear
perturbations are of the form:
\begin{eqnarray}
\label{eqeta}
\eta_1(x) & = & \tilde{\eta_1} e^{i \alpha x} \\
\label{eqc}
\mathbf{u}(x,y) & = & \mathbf{c} e^{\mp |\alpha| y} e^{i \alpha x}
\end{eqnarray}
where (\ref{eqc}) arises from the separable solution of the Laplace
operator in half planes and the sign is chosen appropriately depending
on whether the corresponding variable is in the upper or lower domain,
respectively. Again following \cite{roger}, if the forms (\ref{eqeta})
and (\ref{eqc}) are put in to (\ref{eq34}) the first order term
results in the following algebraic condition:
\begin{equation}
\label{eqwell1}
\mathbf{GMc}+\hat{\eta}_1\mathbf{G}\mathbf{\tilde{U}}_{0\mathbf{n}}=\mathbf{f}
\end{equation}
where
\[
\mathbf{M}=\left( \begin{array}{cccc} -|\alpha| &0 &0 &0 \\ 0
&|\alpha| &0 &0 \\ 0 &0 &-|\alpha| &0\\0 &0 &0 &|\alpha|
\\ 1 &0
&0 &0 \\ 0  &1 &0 &0 \\ 0 &0 &1 &0\\0 &0 &0 &1
\end{array}
\right)
\]
and
\[
\mathbf{\tilde{U}}_{0\mathbf{n}} =\left(\begin{array}{c}
\mathbf{u}_{0yy} \\ \mathbf{u}_{0y} \end{array}\right)
\]
evaluated at $y=0$. This can be written as:
\[
\mathbf{\tilde{U}}_{0\mathbf{n}}= \mathcal{S} \tilde{\mathbf{U_0}}
\]
where $\tilde{\mathbf{U_0}}$ is the vector from the base solution
(\ref{gub}) and $\mathcal{S}$ is a shift matrix that moves rows 1-4 to
rows 5-8 and then zeros the first four rows.

Considering (\ref{eqwell1}), the well-posedness criteria can be
reduced to the existence and uniqueness of $(\mathbf{c},
\hat{\eta}_1)$ for every $\mathbf{f}$ and $\alpha$. Consider the
matrix $\mathbf{GM}$. It is a $5\times 4$ matrix with rank 4. We
denote a vector in its left nullspace by
$\mathbf{w}=\mathcal{N}[(\mathbf{GM})^T]$. Thus to ensure the
well-posedness of the problem, we require:
\begin{equation}
\label{eqwell2}
\mathbf{w^T}\mathbf{G}\mathbf{\tilde{U}}_{0\mathbf{n}}\neq 0
\end{equation}
This is the algebraic condition for linear well-posedness for the
general class of steady free boundary value problems we have
considered.

\subsubsection{Residual velocity choice}

In our numerical algorithm, we solve a fixed boundary problem each
iteration with some of the boundary conditions in $\mathbf{G}$
satisfied, while the residual in the remaining condition is used as
the normal velocity to evolve the interface.  With the well-posedness
criteria of the problems established above, we can discuss the
stability of the numerical methods for the free boundary problem. We
consider residual velocity methods as discussed in the introduction.
The only time dependence is in the interface position $y = \eta(x,t)$
which in the linearized setting above reduces to
\[
y = \tilde{\eta} e^{i \alpha x} e^{\lambda t}
\]
The stiffness $\lambda$ is a function of the choice of residual to use
as velocity and also the wave number $\alpha$. If for a given residual
velocity choice, $ Re(\lambda)<0$ for all $\alpha$, the resulting
method will converge to the steady interface solution. We show how in
addition residual velocities can be chosen in some cases so that
$\lambda$ is bounded independently of $\alpha$, which suggests that
time steps for the resulting method can be chosen independently of the
spatial discretization. This is verified in the example computation in
the next section.

Again, following \cite{roger} the linear analysis leads to
\begin{equation}
\label{eqlambda}
\lambda=\frac{\mathbf{w^T
G\mathbf{\tilde{U}}_{0\mathbf{n}}}}{\mathbf{w^T v}}
\end{equation}
where $\mathbf{v}$ denotes the linear combination of residuals used as
the normal velocity (the orthogonal combinations are used as interface
conditions at each time) and as before $\mathbf{w^T}$ is the left
nullspace of $\mathbf{GM}$.

Note that in the numerator of (\ref{eqlambda}) is the well-posedness
term (\ref{eqwell2}). We assume we are only trying to compute
solutions to well-posed problems, so can assume this term does not
change sign. Thus, the stability of the scheme depends only on the
denominator which depends on a straight-forward way on the velocity
choice $\mathbf{v}$.

\subsubsection{Discussion}

The well-posedness condition (\ref{eqwell2}) and stability condition
(\ref{eqlambda}) are powerful analytic tools. The influence of the
boundary conditions (through $\mathbf{G}$), the elliptic problem
(through $\mathbf{w}$ coming from $\mathbf{M}$), and boundary and
far-field data (through $\mathbf{\tilde{U}}_{0\mathbf{n}}$) are
encoded in these scalar relationships. However, for problems of size
``2+2'' and larger, it becomes difficult to make general statements
about a particular problem. In general, for ``2+2'' problems,
$\mathbf{w}$ has elements that are fourth order polynomials in
$|\alpha|$. Thus, (\ref{eqwell2}) will in general also be a fourth
order polynomials in $|\alpha|$, with coefficients depending on
boundary data and global fluxes. We discuss below situations of
physical interest (including the main example from the last section)
where significant simplification occurs in these expressions.

\subsection{Application to physical example}

Let's go back to the concrete example in section \ref{sphys} where
$\mathbf{G}$ is given in (\ref{eq28}). The nullspace of
$\mathbf{G}$, $\mathcal{N}(\mathbf{G})$, is:
\begin{equation}
\mathcal{N}(\mathbf{G})^T=\left(
\begin{array}{ccccccccc}
  &\frac{K^-}{K^+} &1 &1 &0  &0 &0 &0 &0 \\
 &-\frac{1}{K^+} &0  &+1 &+1  &0   &0 &0 &0 \\
 &0 &0  &0  &0   &1   &1 &1 &0
 \end{array} \right)
\end{equation}
Each row of the matrix above corresponds to the effect of one global
flux. We describe the physical interpretations of these global fluxes
below:
\begin{description}
\item[$q_1$:] This corresponds physically to a heat flux downwards
in the lower region. This heat flux is balanced by a heat flux in the
upper region. In addition, the temperature gradient in the lower
region drives a vapour flux there that is matched by a pressure driven
gradient flow in the upper region.
\item[$q_2$:] This corresponds physically to a flux of water downwards
in the lower region. This flux is generated by condensation at the two
phase boundary which is provided by vapour flow downwards and heat
flow upwards in the upper region.
\item[$q_3$:] This term represents an additive shift in temperature and
pressure.
\end{description}
The data term $\mathbf{b}$ is zero in this example, so the
well-posedness is determined entirely by the global fluxes.
Specifically, the term $\mathbf{\tilde{U}}_{0\mathbf{n}}$ in
(\ref{eqwell2}) and (\ref{eqlambda}) is given by
\[
\mathbf{\tilde{U}}_{0\mathbf{n}} \in
\mathcal{S}[\mathcal{N}(\mathbf{G})]=\left(
\begin{array}{cccccccc}
 0 &0 &0 &0 &\frac{K^-}{K^+} &1 &1 &0 \\
 0 &0 &0 &0 &-\frac{1}{K^+} &0 &+1 &+1 \\
 0 &0 &0 &0 &0 &0 &0  &0
 \end{array} \right)^T
\]
where by the last term above we mean the span of the column vectors,
which has rank 2. This is our first simplification: that only two of
the global fluxes (the first two described above) enter the stability
condition. Upon reflection, it is clear that shifting the overall
temperature and pressure (the effect of $q_3$) in this linear problem
should not affect the stability of the problem.

We now find the left nullspace $\mathbf{w}$ of
$\mathbf{(GM)}^T$:
\[
\mathbf{w^T}=( \frac{\,K^++K^-+2}{K^++K^-}|\alpha|,\;\frac{-2
K^+}{K^++K^-}|\alpha|,\;-|\alpha|,\; \frac{-2}{K^++K^-},
 1\;)
\]
The first three elements are $O(|\alpha|)$, corresponding to rows in
$\mathbf{G}$ representing pure Dirichlet conditions. This is the
second simplification observed in this example, that the vector
$\mathbf{w^T}$ has relatively simple dependence on $|\alpha|$.

By simple calculation, we get:
\[
\mathbf{w^T}\mathbf{G}\;\mathbf{\tilde{U}}_{0\mathbf{n}}=2\;|\alpha|\;\frac{(K^+-K^-)q_1+(K^++K^-+2)q_2}{K^++K^-}
\neq 0
\]
to be the well-posedness criteria, that is
\begin{equation}
\label{eq311}
(K^+-K^-)q_1+(2 + K^+ + K^-)q_2 \neq 0.
\end{equation}
The simplifications noted above are what lead to this particularly
clear form. We note that in the physical situation described in
\cite{lloyd}, $q_2 < 0$ (liquid moves toward the two phase zone) and
$q_1 \approx 0$ so we expect the well-posedness term in (\ref{eq311})
to be negative.

We consider now the choice of residual velocity (\ref{eqlambda}) for
this case. For instance, if $\mathbf{v}= \mathbf{e}_5$ (the physical
velocity choice), we have:
\[
\lambda =  2\;|\alpha|\;\frac{(K^+-K^-)q_1+(K^++K^-+2)q_2}{K^++K^-}.
\]
This is a stable choice for velocity, considering the predicted sign
of the term (\ref{eq311}) above. However, when we consider $\alpha \to
\infty$ (finer discretizations), we get a very stiff problem and the
time step we take should be correspondingly small. Therefore, better
numerical properties are obtained when we choose $\mathbf{v}$ such
that $\lambda$ is independent of wave number. This goal can be
achieved by choosing residuals of Dirichlet conditions. In this case,
we can choose $\mathbf{v}^T= \mathbf{e}_1$ (corresponding to the
saturation condition), then
\[
\lambda = 2\;\frac{(K^+-K^-)q_1+(K^++K^-+2)q_2}{K^++K^- +2}
\]
and we have an optimal velocity to evolve the interface. Numerical
validation of the superior performance of this residual velocity over
the physical one is given in Section \ref{snum} below.

A straight forward calculation shows that for all problems of class
$\tilde{C}$ identified at the end of section \ref{sclass} similar
results are obtained. That is, the vector $\mathbf{w}$ in
(\ref{eqwell2}) and (\ref{eqlambda}) has terms of order $\alpha$ in
entries corresponding to pure Dirichlet conditions and constant in
entries corresponding to pure Neumann conditions. Also, only the two
global fluxes corresponding to derivative quantities enter these
expressions. Similar simplifications are found in other cases when the
interface conditions separate into pure Dirichlet and pure Neumann
conditions.

\subsection{Extension to the 3-dimensional problem}

The linear analysis is similar except that we assume an initial flat
interface to be $z=0$, and a small perturbation of driving function
is given:
\begin{equation}
\label{eq312}
\mathbf{G\tilde{U}}=\mathbf{b}+\epsilon \mathbf{f} e^{i\alpha x}
e^{i\beta y}
\end{equation}
then consequently we have updated solution and interface:
\[
\eta \sim \eta_0+\epsilon \eta_1+ \cdots
\]
\[
\mathbf{u} \sim \mathbf{u}_0+\epsilon \mathbf{u}_1+\cdots
\]
where
\[
\eta_1(x,y)=\hat{\eta}_1 e^{i\alpha x} e^{i\beta y}
\]
Then the separable solutions of Laplacian equation $\Delta
\mathbf{u}=0$ are in form of:
\[
\mathbf{u}=e^{i\alpha x}e^{i\beta y} \mathbf{\hat{u}}
\]
where
\[
\mathbf{\hat{u}}=\mathbf{c}e^{\mp \sqrt{\alpha^2+\beta^2}\;z}
\]
and $\mathbf{c}=(A^+,\;A^-,\;B^+,\;B^-)^T$ is any arbitrary constant
vector determined by the fixed boundary conditions. So we have:
\begin{equation}
\label{eq313}
\mathbf{u}_1=\mathbf{c}e^{\mp \sqrt{\alpha^2+\beta^2}\;z}e^{i\alpha
x} e^{i\beta y}
\end{equation}

If we expand the solution along the new interface, we have:
\begin{equation}
\label{eq314}
\mathbf{u}(x,y,\eta(x,y))=\mathbf{u}_0(x,y,0)+\epsilon \eta_1
\mathbf{u}_{0z}(x,y,0)+\epsilon \mathbf{u}_1(x,y,0)+O(\epsilon^2)
\end{equation}
and since the leading order term of interface is $\eta_0=0$, we can
consider roughly $\mathbf{u}_\mathbf{n}=\mathbf{u}_z$.
\begin{equation}
\label{eq315}
\mathbf{u}_\mathbf{n}(x,y,\eta(x))=\mathbf{u}_{0z}(x,y,0)+\epsilon[\mathbf{u}_{1z}(x,y,0)+\eta_1\mathbf{u}_{0zz}(x,y,0)]+O(\epsilon^2)
\end{equation}
Substitute those expansions into (\ref{eq314}), equating coefficients of
different orders, we have:
\begin{equation}
\label{eq316}
[\mathbf{G\tilde{U}}_1+\hat{\eta}_1 e^{i\alpha x} e^{i\beta
y}\mathbf{\tilde{U}}_{0\mathbf{n}}]=\mathbf{f}\,e^{i\alpha x}
e^{i\beta y}
\end{equation}
where
$\mathbf{\tilde{U}}_{0\mathbf{n}}=(\mathbf{u}_{0zz}^T,\mathbf{u}_{0z}^T)^T$
is dependent on $\mathbf{u}_0$, and
$\mathbf{\tilde{U}}_1=(\mathbf{u}_{1z}^T, \mathbf{u}_1^T)^T$ with
$\mathbf{u}_1$ in (3.15).

The base solution $\mathbf{u}_0$ in 3-D retains the 
same form in the 2-D case.
We need the solutions independent of coordinate
variables $x$ and $y$, which lead to
$\mathbf{u}_0=z\mathbf{r}^1+\mathbf{r}^0$, and consequently
\[(\mathbf{r}^1,\;\mathbf{r}^0)^T \in \mathcal{N}(\mathbf{G})\]
we always have
$\mathbf{\tilde{U}}_{0\mathbf{n}}=(0,\mathbf{r}_1)^T$. Hence:
\begin{equation}
\label{eq317}
\mathbf{\tilde{U}}_{0\mathbf{n}} \in
\mathcal{S}[\mathcal{N}(\mathbf{G})]=\left(\; \mathbf{0}_{3\times
4},\;\; \mathcal{N}(\mathbf{G})_{l4}\; \right)^T
\end{equation}
where $\mathcal{N}(\mathbf{G})_{l4}$ represents the left four
columns of $\mathcal{N}(\mathbf{G})$.

After we figure out $\mathbf{\tilde{U}}_{0\mathbf{n}}$, the next
step is to find $\mathbf{\tilde{U}}_1$. We reserve the notation of
matrix $\mathbf{M}$ to provide us a more clear structure of
$\mathbf{\tilde{U}}_1$.
\[
\mathbf{\tilde{U}}_1(x,y,0)=(\mathbf{u}_{1z},\mathbf{u}_{1})=
e^{i\alpha x}e^{i\beta y}\mathbf{Mc}
\]
and
\[
\mathbf{M}=\left( \begin{array}{c} \mp
\sqrt{\alpha^2+\beta^2}I_{4}\\I_{4}
\end{array}\right)\]
where $I_4$ is a $4\times 4$ identity matrix. Thus the equation
(3.18) becomes
\begin{equation}
\label{eq318}
\mathbf{GMc}+\hat{\eta}_1\mathbf{G}\mathbf{\tilde{U}}_{0\mathbf{n}}=\mathbf{f}
\end{equation}
Denoting its left nullspace by
$\mathbf{w}=\mathcal{N}[(\mathbf{GM})^T]$, we still have the same
criteria for well-posedness:
\begin{equation}
\label{eq319}
\mathbf{w^T}\mathbf{G}\mathbf{\tilde{U}}_{0\mathbf{n}}\neq 0
\end{equation}
Applying this to our physical example in 3-D, we have:
\[
\mathbf{w^T}\mathbf{G}\;\mathbf{\tilde{U}}_{0\mathbf{n}}=2\;\sqrt{\alpha^2+\beta^2}\;\frac{(K^+-K^-)q_1+(K^++K^-+2)q_2}{K^++K^-}
\neq 0
\]
This provides the same information as the 2-D problem. A study of the 
properties of interface velocities also gives results analogous to the 
2-D case.


\section{A 2-D computational example}
\label{snum}

In this section, we show how the residual velocities work on the
vector Laplacian problem by computing the linear example of section
\ref{sphys}. In addition to the interface conditions, the following
fixed boundary conditions are given: On the bottom, we assume the
medium is saturated with liquid water, thus $s = 1$, and give a
constant reference temperature
\[
T^-(x, 0) = T_0.
\]
On the top boundary, we make the boundary impenetrable to vapour
\[
\frac{\partial P}{\partial \mathbf{n}} = 0
\]
and take the heat flux to be given:
\[
K^+ \frac{\;\partial T^+}{\partial \mathbf{n}} = f_1 (x).
\]
We take $f_1(x) = 2 + \sin(x)/2$ and $T_0 = 10$ as our numerical
example. These conditions lead to a global flux $q_2 <0$ (water flux 
upwards towards the interface) at each 
interface location. In the computations below, we take $K^+=K_-$ 
so the well-posedness condition (\ref{eq311}) is satisfied with the 
sign that makes our velocity sign choices stable. 

\subsection{Exact solution}

In order to get an exact solution that can be compared with
the numerical approximation, we start by building two piecewise
functions
\[U = \begin{cases}  &K^+ T^+
\qquad\mbox{in}\;\;D^+ \\ & K^-T^- - s \quad\mbox{in}\;\;D^-
\end{cases} \qquad  V = \begin{cases}  &P
\qquad\mbox{in}\;\; D^+ \\&T^- + s   \quad\mbox{in}\;\;D^-
\end{cases}  \]
The two Laplacian equations $\Delta U = 0$ and $\Delta V = 0$ hold
in the whole domain discarding the interface. The continuity gives
the Neumann conditions across the interface:
\[
K^+ T^+_\mathbf{n} - K^- T^-_\mathbf{n} = -s_\mathbf{n} \qquad
\mbox{and} \quad P_\mathbf{n} = T^-_\mathbf{n}+ s_\mathbf{n}
\]
and because $s = 0$ at the interface, the Dirichlet conditions across
the interface are:
\[
K^+ T^+ = K^- T^- \qquad \mbox{and} \quad P = T^-
\]
Note that only when $K^+ = K^-$, the problem is equivalent to our
original one, so we set $K^+ = K^-$ and try to find the location of
the interface $s=0$. Solving the problem, we have:
\[
\begin{aligned}
T^+ &= 9 + 2y + \frac{\sin(x)}{2} \frac{\sinh
(y)}{\cosh(L)} \qquad  P = 11\\
T^- &= 10 + y + \frac{\sin(x)}{4} \frac{\sinh (y)}{\cosh(L)} \qquad
s = 1 - y - \frac{\sin(x)}{4}  \frac{\sinh (y)}{\cosh(L)}
\end{aligned}
\]
and the location of $s = 0$ is given implicitly by
\[
4(1-y)\cosh(L) - \sin(x)\sinh(y) = 0
\]
By setting $x$ equally distributed in the interval $[0, 2\pi]$, we
can solve this transcendental equation for $y$, and thus find the
coordinates of points on the interface to arbitrarily high 
precision. 

\subsection{Finite difference method in mapped rectangular domain}

To handle the two irregular subdomains separated by curved
interface, we map them into two rectangular domains, so that finite
difference method can be naturally employed. A finite
element method would be a good alternative for more irregular
domains. A preliminary finite element implementation was done in
\cite{roger} which required considerable grid refinement near the
interface. The finite difference implementation given here
demonstrates that this is not necessary for accurate solutions.

Suppose $y = h(x)$ is the known interface, we will do the following
mapping:
\[ \begin{aligned}
& D^+:\;\;(x, y)\;\to (x_1 = x,\;\; y_1 = 1 + (y - h(x))/(L-h(x)) ) \\
& D^-:\;\;(x, y)\;\to (x_2 = x,\;\; y_2 = y/h(x) )
\end{aligned}
\]
The new domain is $[0, 2\pi] \times [0, 2]$, and the interface $y =
h(x)$ has been mapped into $y_1=y_2=1$. Then by simple calculation,
we have:
\begin{eqnarray*} \frac{\partial T}{\partial x} &=&
\frac{\partial T}{\partial x_2} \frac{\partial x_2}{\partial x} +
\frac{\partial T}{\partial y_2}\frac{\partial y_2}{\partial x} =
\frac{\partial T}{\partial x_2} - \frac{h'(x_2)y_2}{h(x_2)}
\frac{\partial T}{\partial y_2}\\
 \frac{\partial T}{\partial y} &=&
\frac{\partial T}{\partial x_2}\frac{\partial x_2}{\partial
y}+\frac{\partial T}{\partial y_2}\frac{\partial y_2}{\partial y} =
\frac{1}{h(x_2)}\frac{\partial T}{\partial y_2}
\end{eqnarray*}
Note that $\frac{\partial x_2}{\partial y} = 0, \frac{\partial
y_2}{\partial x} \neq 0$. Rewrite the Laplacian equations in the
mapped lower subdomain: \small
\begin{equation}
\begin{aligned}
\Delta T^- &= \frac{\partial^2 T}{\partial y^2}+\frac{\partial^2
T}{\partial x^2} = \frac{1}{h(x_2)^2}\frac{\partial^2 T}{\partial
y_2^2} +\frac{\partial x_2}{\partial x} \frac{\partial}{\partial
x_2} \left(\frac{\partial T}{\partial x}\right) -
\frac{h'(x_2)y_2}{h(x_2)}\frac{\partial y_2 }{\partial x}
\frac{\partial}{\partial y_2} \left(\frac{\partial T}{\partial
x}\right)
\\
&= \frac{1}{h(x_2)^2}\frac{\partial^2 T}{\partial y_2^2}
+\frac{\partial x_2}{\partial x} \frac{\partial}{\partial x_2}
\left( \frac{\partial T}{\partial x_2} - \frac{h'(x_2)y_2}{h(x_2)}
\frac{\partial T}{\partial y_2} \right) \\
&- \frac{h'(x_2)y_2}{h(x_2)}\frac{\partial y_2 }{\partial x}
\frac{\partial}{\partial y_2} \left( \frac{\partial T}{\partial x_2}
- \frac{h'(x_2)y_2}{h(x_2)} \frac{\partial T}{\partial y_2} \right)
\\ &= \frac{1+h'(x_2)^2 y_2^2}{h(x_2)^2} \frac{\partial^2
T^-}{\partial y_2^2} + \frac{\partial^2 T^-}{\partial x_2^2} -
2\frac{y_2 h'(x_2)}{h(x_2)}\frac{\partial^2 T^-}{\partial x_2
\partial y_2}\\ &+ \frac{2 h'(x_2)^2 - h''(x_2) h(x_2)}{h(x_2)^2} y_2 \frac{\partial T^-}{\partial y_2}
\end{aligned}
\end{equation}
Similar calculation gives the Laplacian equation in the mapped upper
subdomain:
\begin{equation}
\begin{aligned}
\Delta T^+ &= \frac{1+h'(x_1)^2 (2-y_1)^2}{(L-h(x_1))^2}
\frac{\partial^2 T^+}{\partial y_1^2} + \frac{\partial^2
T^+}{\partial x_1^2} - 2\frac{(2-y_1)
h'(x_1)}{(L-h(x_1))}\frac{\partial^2 T^+}{\partial x_1
\partial y_1} \\&- \frac{2 h'(x_1)^2 + h''(x_1)(L-h(x_1))}{(L-h(x_1))^2} (2-y_1) \frac{\partial T^+}{\partial y_1}
\end{aligned}
\end{equation}
\normalsize For the pressure $P$ and saturation $s$, we get similar
formulas. The Neumann interface and fixed boundary conditions should
be reformulated as well in the new coordinates. 
For the
Neumann (physical) velocity computation, we choose the
residual of mass conservation condition $P_\mathbf{n} -
T^-_\mathbf{n}- s_\mathbf{n} = 0$ to move the interface.
Alternatively, the Dirichlet condition $s = 0$ is picked as the
Dirichlet residual velocity.

In order to make our numerical approximation second order accurate,
we use central difference scheme all through the discretization.
With $N\times N$ to be the resolution on either side, the resultant
coefficient matrix is of order $4N^2 \times 4N^2$.

\subsection{The numerical results}

We start from $N=10$, and increase it to get higher accuracy. $L=2$
is fixed all through the computation. To get the results in Table
\ref{t1}, we set the timestep $\Delta t=0.02$ for the Neumann
(physical) velocity, and compute till $T = 24$, which is long enough
to reach the steady state. The errors \textbf{errInf, errT, errS}
are respectively the errors of interface, temperature and saturation
in maximum norm. For the Dirichlet (artificial) residual velocity,
we take $\Delta t=0.2$, and iterate till $T=24$, the errors are
shown in Table \ref{t2}. Note that at (almost) steady state, 
the errors of the two methods are (almost) identical since the 
discretization of the steady-state solution is identical. 

\begin{table}
\begin{center}
\begin{tabular}[b]{|c||c|c|c|c|c|c|}
\hline Grids   &errInf  &ratio &errT &ratio    &errS &ratio \\
\hline\hline
 10$\times$10 &0.0031  & ***  &0.0075 & ***  &2.8773e-5 & ***  \\ \hline
 20$\times$20 &8.0341e-4 &3.859 &0.0019 &3.947    &6.4194e-6 &4.48  \\ \hline
 40$\times$40 &2.0019e-4 &4.013 &4.8504e-4 &3.917   &1.6205e-6 &3.961  \\ \hline
\end{tabular}
\caption{\label{t1} Convergence test of Neumann residual velocity}
\end{center}
\end{table}

\begin{table}
\begin{center}
\begin{tabular}{|c||c|c|c|c|c|c|} \hline
Grids        &errInf  &ratio   &errT &ratio    &errS &ratio \\
\hline\hline
10$\times$10  &0.0032 & *** &0.0075 & ***  &2.9558e-5 & *** \\
\hline 20$\times$20  &8.0341e-4 &3.983 &0.0019 &3.947  &6.4127e-6 &4.601 \\
\hline 40$\times$40 &2.0019e-4 &4.013 &4.8504e-4 &3.917 &1.6201e-6
&3.960
\\ \hline
\end{tabular}
\caption{\label{t2} Convergence test of Dirichlet residual velocity}
\end{center}
\end{table}

Obviously, both methods have second order accuracy implied by the
convergent rate. Now we try to demonstrate the advantage of the
Dirichlet (artificial) residual velocity by setting the timestep
$\Delta t$ a free parameter. As we know, since we have to use explicit
scheme for $\mathbf{x}_{\mathbf{n}} = R(\mathbf{x})$ to evolve the
interface, the timestep has to be reduced when refining the grid to
make it within the stability region. In Table \ref{t3}, we observe
this property in Neumann residual velocity, but not in Dirichlet
one. This means, even when we have very fine grids, the timestep can
be kept reasonably large. Since we know in each time iteration, the
computational complexity is almost the same, allowing larger
timesteps can greatly reduce the computational cost.

In Table \ref{t3}, $\textbf{VN}$ and $\textbf{VD}$ represent
respectively the Neumann and Dirichlet residual velocities. It's not
hard to see for Neumann residual velocity, $\Delta t = 0.12$ is
appropriate for $10 \times 10$ grids to get convergence, and when we
refine to $20\times20$, $\Delta t = 0.06$ is correspondingly
decreased, and for $40\times40$, it becomes $\Delta t =
0.03$. However, we can use $\Delta t = 0.2$ for all three resolutions
if the Dirichlet residual velocity is chosen. Considering the
appropriate timesteps $\Delta t = 0.03$ and $\Delta t = 0.2$ when the
resolution is $40\times40$, Dirichlet residual velocity is obviously
better than the physical Neumann one.

\begin{table}
\begin{center}
\begin{tabular}{|c||c|c|c|c|c|} \hline
Grids       &timestep     &RV     &errInf     &errT    &errS  \\
\hline\hline
10$\times$10 &0.20        &VN  &0.9213  & 0.8734   &0.1029  \\
\cline{2-6}  &0.15        &VN  &0.0453  & 0.0863   &0.0122  \\
\cline{2-6}  &0.12        &VN  &0.0031  & 0.0075   &2.8773e-5 \\
\hline
10$\times$10 
&0.20        &VD  &0.0032  & 0.0075   &2.9558e-5 \\
\hline\hline
20$\times$20  &0.12        &VN  &12.0340  &6.2583   &0.8346  \\
\cline{2-6}  &0.08        &VN  &73.4589  &4.9261   &1.3080 \\
\cline{2-6}  &0.06        &VN  &8.0341e-4  &0.0019   &6.4194e-6 \\
\hline
20$\times$20 
&0.20        &VD  &8.0341e-4  &0.0019   &6.4127e-6 \\
\hline\hline
40$\times$40  
&0.06        &VN  &9.6403 &24.4088 &0.9594 \\
\cline{2-6}  &0.04       &VN  &3.8939 &3.5538 &0.8762 \\
\cline{2-6}  &0.03       &VN  &2.0019e-4 &4.8504e-4 &1.6205e-6 \\
\hline
40$\times$40 
&0.20        &VD  &2.0019e-4 &4.8504e-4 &1.6201e-6 \\
\hline
\end{tabular}
\caption{\label{t3} Performance of different residual velocities}
\end{center}
\end{table}

\section{Conclusions and discussion}

Local, linear well-posedness criteria were developed for a general
class of ``2+2'' vector Laplacian problems. Much simplification can be
obtained in certain cases of physical interest. This is demonstrated
in a simplified model of two-phase flow with phase change in porous
media. The theory also allows the investigation of the stability and
numerical properties of residual velocities to compute the steady
state. It was possible to identify a superior velocity choice for the
physical example and other problems in its general class. The theory
was verified computationally using a finite difference method with
mapped coordinates. In general, this approach offers a way to improve
the performance of certain steady-state interface computations with
very little effort. Ongoing work suggests that these ideas can be
extended to interface problems involving more general second order
elliptic systems, such as problems from viscous, incompressible fluid
flow. Extension to higher order problems with biharmonic operators are
also possible. Efficient implementation of the method in a finite
element framework is still an open problem.

\section*{Acknowledgement}

This work is supported by NSERC and MITACS. We also acknowledge the
valuable suggestions from Lloyd Bridge and Roger Donaldson.


\begin{thebibliography}{99}
\bibitem{al} N. Al-Rawahi and G. Tryggvason (2004) 
\textit{Numerical simulation of dendritic solidification with 
convection: three-dimensional flow}, 
J. Comput. Phys. 194:677-696.

\bibitem{alt} H.W.Alt and L.A.Caffarelli (1981) \textit{Existence and
regularity for a minimum problem with free boundary}, Journal fur die
Reine und Angewandte Mathematik 325:105-144.

\bibitem{lloyd} L.Bridge, R.Bradean, M.J.Ward and B.R.Wetton (2003)
\textit{The analysis of a two-phase zone with condensation in a porous
medium}, Journal of Engineering Mathematics, 45:247-268.

\bibitem{caflisch} R. Caflisch, F. Gibou, R. Fedkiw and S. Osher
(2003) \textit{A level set approach for the numerical simulation of
dendritic growth}, J. Sci. Comput. 19:183-199.

\bibitem{droplet1} T. Cheng, P. Minev, and K. Nandakumar (2004)
\textit{A projection scheme for incompressible multiphase
flow using adaptive Eulerian grids}, Int. J. Numer. Meth. Fluids, 
45: 1-19.

\bibitem{crank} J.Crank (1984) \textit{Free and moving boundary
problems}, Clarendon, Oxford.

\bibitem{roger} R.Donaldson and B.R.Wetton (2004) \textit{Solving
steady interface problem using residual velocities}, accepted in the
IMA Journal of Applied Mathematics, preprint available at 
www.math.ubc.ca/$\sim$wetton/. 

\bibitem{direct} J.Glimm, J.W.Grove and Y. Zhang (2002) 
\textit{Interface tracking for
axisymmetric flows}, Journal of Scientific Computation, 24(1):
208-236

\bibitem{shape2} J.Haslinger, T.Kozubeck, K.Kinisch and G.Peichl (2003)
\textit{Shape optimization and fictitious domain approach for
solving free boundary problems of bernoulli type}, Computational
Optimization and Applications, 26:231-251.

\bibitem{modified} Z.Li (1997) \textit{Immersed interface method for moving interface
problems}, Numerical Algorithms 14:269-293

\bibitem{level} S.Osher and J.A.Sethian (1988) \textit{Fronts propagating with curvature-dependent
speed: Algorithms based on Hamilton-Jacobian formulations}, Journal
of Computational Physics, 100:209-228

\bibitem{bubble} F.Raymond and J.M. Rosant (2000)
\textit{A numerical and experimental study of the terminal 
velocity and shape of bubbles in viscous liquids},
 Chem. Eng. Sci. 55:943-955.

\bibitem{shape1} J.Sokolowski and J.P.Zelesio (1992) \textit{Introduction
to shape optimization, shape sensitivity analysis}, Springer-Verlag,
New York.

\bibitem{droplet2} M. Sussman, P. Smereka and S. Osher (1994) 
\textit{A level set approach for computing solutions to incompressible
2-phase flow}, J. Comput. Phys. 114:146-159.

\bibitem{value} T. Tiihonen (1997) \textit{Shape optimization and 
trial methods for free boundary problems}, RAIRO, Math. Model.
Numer. Anal. 31:805-825.

\bibitem{tracking1} S.O.Unverdi and G.Tryggvason (1992) \textit{A front tracking method
for viscous incompressible, multi-fluid flows}, Journal of
Computatioanal Physics, 100:25-37

\bibitem{trackcap} Y.C.Tai, S.Noelle, J.M.N.T.Gray, and K.Hutter
(2002) \textit{Shock-Capturing and front tracking methods for
Granular Avalanches}, Journal of Computational Physics 175,269-301

\bibitem{tracking} D.E.Womble (1989) \textit{A front-tracking
method for multiphase free boundary problems}, Journal on Numerical
Analysis, 26(2):380-396

\end{thebibliography}
\end {document}